%% file: Moduli_Spaces.tex
\documentclass[11pt,a4paper,twoside]{amsart}
\usepackage{amsmath,amssymb,amsfonts,amsthm,mathrsfs}
\usepackage{times,hyperref,color}

\let\mathcal\mathscr

\newtheorem{The}{Theorem}[section]

\newtheorem{Theorem}{Theorem}[section]
\newtheorem{Proposition}[The]{Proposition}
\newtheorem{Lemma}[The]{Lemma}

\theoremstyle{definition}

\newtheorem{Remark}[The]{Remark}



\begin{document}

%\large

\title{
Moduli spaces of model real submanifolds:
\\
Two alternative approaches }

\author{Masoud Sabzevari}
\address{Department of Pure Mathematics,
University of Shahrekord, 88186-34141 Shahrekord, IRAN}
\email{sabzevari@math.iut.ac.ir}

\date{\number\year-\number\month-\number\day}

\maketitle

\begin{abstract}
Instead of the invariant theory approach employed by Beloshaoka and Mamai for constructing the moduli spaces of Beloshapka's universal CR-models, we consider two alternative approaches borrowed from the theories of equivalence problem and Lie symmetries, each of them having its own advantages. Also the moduli space $\mathcal M(1,4)$ associated to the class of universal CR-models of CR-dimension 1 and codimension 4 is computed by means of the presented methods.
\end{abstract}

\pagestyle{headings} \markright{Moduli spaces of model real submanifolds: two alternative approaches }

\section{Introduction}

After 1907 that Henri Poincar\'{e} \cite{Poincare} studied real
submanifolds in the specific complex space $\mathbb C^2$ according to
the associated {\it model surface}, namely the {\it Heisenberg sphere}, the general issue of
investigating real submanifolds in arbitrary complex spaces according
to their associated models gained its increasing interest
\cite{Beloshapka2004,Chern-Moser, 5-cubic, Merker-Sabzevari-CEJM}. In
this approach, many questions about automorphism groups,
classification, invariants and others, concerned the (holomorphic)
transformations of real submanifolds in a certain complex space can
be reduced to similar problems about the associated models.

Developing along the Poincar\'{e}'s approach, Chern and Moser in their famous paper \cite{Chern-Moser} investigated it in the case of real hypersurfaces of
arbitrary dimensions. But as the most general investigation\,\,---\,\,to the best of the author's knowledge\,\,---\,\,
Valerii Beloshapka has
studied extensively the subject of model surfaces in arbitrary dimensions and codimensions and found some
considerable results in this respect. Specifically in 2004,
he introduced and established in \cite{Beloshapka2004} the structure of
some nondegenerate models associated (uniquely) to totally
nondegenerate germs of arbitrary Cauchy-Riemann (CR for short)
dimensions and codimensions. Two such germs are holomorphically equivalent whenever their associated models are equivalent.
He also developed a machinery of the
construction of these models ({\it see also} \cite{Mamai}). Beloshapka called each model $M\subset\mathbb C^{n+k}$
of certain CR-dimension and codimension $n$ and $k$ by a {\sl universal} CR-manifold of the {\sl type}
$(n,k)$ which is homogeneous and enjoys several {\it nice}
properties (\cite[Theorem 14]{Beloshapka2004}) that
exhibit the significance of the models.

Amongst all universal CR-models of various types $(n,k)$, the class of those of CR-dimension $n=1$ has gained more considerations in the literature
 \cite{Beloshapka1997, BES, Merker-Sabzevari-M3-C2, Mamai, Mamai2009, 5-cubic-geometry,  Merker-Sabzevari-CEJM,  5-cubic, Shananina2000}. In particular, the Cartan geometries of the types $(1,1), (1,2)$ and $(1,3)$ are studied in \cite{Merker-Sabzevari-CEJM, BES, 5-cubic-geometry} and the equivalence problems of the totally nondegenerate CR-manifolds associated to the types $(1,1)$ and $(1,3)$ are solved in \cite{Merker-Sabzevari-M3-C2, 5-cubic}.

An important fact about the already mentioned CR-models of the types $(1,k)$ is that in contrary to the cases $k=1,2,3$, almost all the next types $(1,4), (1,5), \ldots$ do not admit a single unique universal CR-model. More precisely, the associated defining equations of theses generic CR-models depend also to some appearing {\it parameters} for which their different values give different CR-models of a fixed type. Therefore for $k\geq 4$, only it is not of interest the holomorphic equivalences between totally nondegenerate CR-manifolds of the same type $(1,k)$ ({\it cf.} \cite{Merker-Sabzevari-CEJM, 5-cubic}) but also it will be interesting to ask under which circumstances on the appearing parameters two universal CR-models of a fixed type are biholomorphically equivalent. For the first time, Beloshapka introduced this question in \cite{Beloshapka} and called the family of all equivalence CR-models of a fixed type $(n,k)$ by the {\sl moduli space} of that type, denoted by $\mathcal M(n,k)$. Subsequently Mamai in \cite{Mamai}, developed this concept by computing the invariants of the moduli spaces $ \mathcal M(1,k)$ for $k=8, \ldots, 13$, where the spaces $\mathcal M(1,4), \mathcal M(1,5)$ and $\mathcal M(1,7)$ were considered by Beloshapka, himself, in \cite{Beloshapka}. The nontrivial properties of the topological structure of the moduli spaces enables one to introduce a series of some new CR-characteristic classes.

As is shown in \cite[Corollary 7]{Beloshapka2004}, two model surfaces are holomorphically equivalent if and only if they are {\it linearly} equivalent. This significant fact has a key role in \cite{Beloshapka} and \cite{Mamai} to establish the method of constructing the desired moduli spaces. Indeed by taking this fact into account, Beloshapka and Mamai employed some powerful techniques from linear algebraic groups and invariant theory (\cite{invariant-theory}), based on the Hilbert basis theorem and rational quotients, for computing the invariants of the desired moduli spaces. We refer the reader to these two papers for more details of this method.

The goal of this paper is to set up two {\it alternative approaches}\,\,---\,\,instead of that from invariant theory\,\,---\,\,for computing the desired invariants of the mentioned moduli spaces, namely the approach borrowed from the theory of equivalence problem and the approach borrowed from the theory of Lie symmetries. Our original impulse for this investigation came from the long memoir \cite{5-cubic}. Each of the three approaches from invariant theory, Cartan equivalence problem, and Lie symmetry has its own advantages.

The first approach comes from the theory of equivalence problem initiated for the first time by \'{E}lie Cartan \cite{Cartan-1932}. The algorithm  devised by Cartan and subsequently developed by Chern and Moser in \cite{Chern-Moser} is a powerful method to solve the biholomorphic equivalence problem between nondegenerate CR-manifolds. The major advantages of using this approach lies in the facts that:
\begin{itemize}
\item[$\bullet$] In this approach we do not need to know any considerable fact or feature of the under consideration CR-models.
\item[$\bullet$] It is more systematic and manageable by computer softwares.
\end{itemize}

The second approach relies on computing the symmetry Lie algebras of the under consideration CR-models, what is called by the Lie algebras of infinitesimal CR-automorphisms in the terminology of CR-geometry. This second approach enjoys the following two advantages:
\begin{itemize}
\item[$\bullet$] By means of the algorithm designed in \cite{SHAM}, one can proceed the most complicated part of the associated computations, namely computing the already mentioned symmetry algebra, by the aid of the known computer softwares like {\sc Maple}.
\item[$\bullet$] Once one computes the associated symmetry Lie algebras, then it is possible to bypass a major part of the Cartan algorithm for finding the final structure equations of the equivalence problem.
\end{itemize}
However, computing the mentioned Lie algebras of infinitesimal CR-automorphisms by means of the classical methods is in fact complicated and time-consuming ({\it cf.} \cite{Beloshapka1997, Mamai2009, Merker-Sabzevari-CEJM, 5-cubic, Shananina2000}). In particular, the size of the computations grows extensively as soon as the number of the CR-dimension or codimension grows even by one. Nevertheless, recently in \cite{SHAM} we have designed a powerful algorithm to compute the desired algebras by employing just some simple techniques of linear algebra and of the modern theory of comprehensive Gr\"obner systems instead of constructing and solving some {\sc pde} systems arising in the classical method.

One should notice that in this paper we do not aim to confirm the Beloshapka and Mamai's results of their papers \cite{Beloshapka, Mamai}, but our aim is to introduce two alternative approaches for computing the invariants of the desired moduli spaces. Then, here we only consider the first appearing moduli space  $\mathcal M(1,4)$.

This paper is organized as follows. In section \ref{preliminaries}, we present required definitions, results and terminology of the theories of Cartan equivalence problem and Lie symmetry. In sections \ref{Cartan}, we employ the Cartan algorithm for constructing the invariance of the moduli space $\mathcal M(1,4)$. Finally in section \ref{approach-2}, we consider the second approach of Lie symmetries for the similar construction. We  observe that the achieved single invariant of the moduli space $\mathcal M(1,4)$ is precisely that computed by Beloshapka and Mamai in \cite{Beloshapka, Mamai}.

\section{Preliminaries}
\label{preliminaries}

\subsection{Cartan equivalence problem}

The main goal in the theory of {\sl equivalence problem}, is to determine whether two given geometric structures can be mapped bijectively onto each other by some diffeomorphism. This problem can be considered in many different contexts, such as equivalences of submanifolds, of differential
equations, of frames, of coframes and of several other geometric structures. Unifying the wide variety of these seemingly different equivalence problems into a potentially universal approach, \'{E}lie Cartan showed that almost all continuous classification questions can indeed be reformulated in terms of specific {\it adapted coframes} ({\em see}~\cite{Olver-1995}).

Seeking an equivalence between coframes usually comprises a certain
initial {\sl ambiguity subgroup} $G\subset {\sf GL}(n)$ related to
the specifc features of the geometry under study. The fundamental
general set up is that, for two given coframes
$\Omega:=\{\omega_1,\ldots,\omega_n\}$ and
$\Omega':=\{\omega'_1,\ldots,\omega'_n\}$ on two certain
$n$-dimensional manifolds $M$ and $M'$, there exists a diffeomorphism
$\Phi:M\longrightarrow M'$ making a geometric equivalence {\em if and
only if} there is a $G$-valued function $g:M\rightarrow G$ such that
$\Phi^\ast(\Omega) = g\cdot\Omega'$.

Cartan's algorithm comprises three interrelated principal aspects:
{\sl absorbtion}, {\sl normalization} and {\sl prolongation}.
Nevertheless, the outcomes of this procedure is often unpredictable.

In brief outline,
starting from:
\begin{equation}
\label{equivalence group} \Omega:=g\cdot\Omega',
\end{equation}
one has to find
the so-called {\sl structure equations} by computing
the exterior differential:
\[
d\Omega
=
dg\wedge\Omega'+g \cdot
d\Omega'.
\]
Inverting~\thetag{\ref{equivalence group}}, one then
has to replace the first
term $dg\wedge\Omega'$ by:
\[
dg\cdot g^{-1}\wedge g\cdot\Omega'=\underbrace{dg\cdot
g^{-1}}_{\omega_{\sf MC}}\wedge\,\Omega,
\]
where $\omega_{\sf MC}$ is the so-called Maurer-Cartan matrix form of the associated structure group $G$:

\[
\aligned \omega_{\sf MC} = &\, \big( (\omega_{\sf MC})_j^i
\big)_{1\leqslant j\leqslant n}^{1\leqslant i\leqslant n}:= \,
\sum_{k=1}^n\, dg_k^i\, \big(g^{-1}\big)_j^k:=\sum_{s=1}^r\,
a^i_{js}\,\alpha^s
\endaligned
\]
that one decomposes according to a basis
$\alpha^1, \dots, \alpha^r$
of left-invariant $1$-forms on $G$, with $r := \rm{dim}_{\mathbb R}\, G$, by means
of certain
constants $a_{js}^i$.
Moreover, one has to express the second term $d\Omega'$ above, which
is a $2$-form, as a combination of the $\omega_j \wedge \omega_k$.
Usually, this step is costful, computationally speaking. Now, the
structure equations are received a form like:
\begin{equation}
\label{equiv-1-general-entries}
d\omega_i
=
\sum_{k=1}^n\,\sum_{s=1}^r\,a_{ks}^i\,\alpha^s\wedge\omega_k
+
\sum_{1\leqslant j<k\leqslant n}\, T_{jk}^i\cdot
\omega_j\wedge\omega_k \ \ \ \ \ \ \ \ \ \ \ \ \
{\scriptstyle{(i\,=\,1\,\cdots\,n)}},
\end{equation}
where, the appearing functions
$T_{jk}^i$, which are called by {\sl torsion coefficients},
usually reveal {\em appropriate invariants} of the geometric structure.

The first two major parts absorbtion and normalization of the Cartan algorithm are based on the following fact:

\begin{Proposition} (see \cite[Proposition 12.6]{5-cubic})
\label{prop-changes}
In the structure equations
\thetag{\ref{equiv-1-general-entries}} one can replace each Maurer-Cartan form $\alpha^s$ and each torsion
coefficient $T^i_{jk}$ with:
\begin{equation}
\label{possible-changes}
\aligned
\alpha^s
&
\longmapsto
\alpha^s+\sum_{j=1}^n\,z^s_j\,\theta^j
\ \ \ \ \ \ \ \ \ \ \ \ \ \ \ \ \ \ \ \ \ \ \ \ \ \ \ \ \
{\scriptstyle{(s\,=\,1\,\cdots\,r)}},
\\
T^i_{jk}&\longmapsto T^i_{jk}
+
\sum_{s=1}^r\,
\big(
a_{js}^i\,z_k^s
-
a_{ks}^i\,z_j^s
\big) \ \ \ \ \ \ \ \ \ \ \
{\scriptstyle{(i\,=\,1\,\cdots\,n\,;\,\,\,
1\,\leqslant\,j\,<\,k\,\leqslant\,n)}},\,
\endaligned
\end{equation}
for some arbitrary functions $z^\bullet_\bullet$
on the base manifold $M$.
\qed
\end{Proposition}

Thus, one does pick the functions-coefficients $z_j^s$ in order
to absorb\,\,---\,\,usually to a constant integer like $0$, $1$ or $i$\,\,---\,\,as many as possible torsion coefficients in the
Maurer-Cartan part, then the remaining, {\sl unabsorbable}, (new,
less numerous) torsion coefficients become {\em true invariants} of
the geometric structure under study. Such absorbtion leads one to determine some certain group parameters in an appropriate way and  {\it hopefully} at the end of several possible normalization-absorbtion loops, one reduces
the structure group $G$ to dimension $0$, getting a so-called {\sl
$\{e\}$-structure}.

But if, as often occurs, it becomes no longer possible after several
absorption-normalization steps to determine a (reduced) set of
remaining group parameters, then one has to add the rest of
Maurer-Cartan forms $\alpha^\bullet$ to the initial lifted coframe
$\Omega$ and to {\sl prolong} the base manifold $M$ as the product
$M^{\sf pr}:=M\times G$.  Surprisingly, Cartan observed that ({\it see} \cite[Proposition 12.1]{Olver-1995}) the
solution of the original equivalence problem can be derived from that
of $M^{\sf pr}$ equipped with a new larger coframe. Then, one has to restart the
procedure {\it ab initio} with such a new prolonged problem. This
initiates the third essential feature of the equivalence algorithm:
the {\sl prolongation}. For a detailed presentation of Cartan's
method, the reader is referred to \cite{Olver-1995,5-cubic}.

\subsection{Infinitesimal CR-automorphisms}
\label{preliminary-infinitesimal}

Consider an arbitrary real analytic generic CR-manifold $M\subset\mathbb C^{n+k}$ of CR-dimension $n$ and codimension $k$, represented in coordinates $(z_1,\ldots,z_n,w_1,\ldots,w_k)$ with $w:=u+iv$ as the graph of some:
\[
\aligned
\Xi_j&:=v_j-\Phi_j(z,\overline z,u)=0, \ \ \ \ \ \ {\scriptstyle (j=1,\ldots,k)},
\endaligned
\]
for some real valued polynomial functions $\Phi_\bullet$.
According to definition, a {\sl (local)
infinitesimal CR-automorphism of $M$}, when understood extrinsically,
is a local holomorphic vector field:
\begin{equation}
\label{XX}
{\sf X}
=
\sum_{i=1}^n\,Z^i(z,w) \frac{\partial}{\partial
z_i} + \sum_{j=1}^k\,W^j(z,w) \frac{\partial}{\partial w_j}
\end{equation}
whose real part ${\rm Re}\, {\sf X} = \frac{ 1}{ 2} ( {\sf X} +
\overline{\sf X} )$ is tangent to $M$, namely $( {\sf X} +
\overline{\sf X} )|_{\Xi_j}\equiv 0$ for each $j=1,\ldots,k$.
The collection of all infinitesimal CR-automorphisms of $M$
constitutes a Lie
algebra which is called the {\sl Lie algebra of infinitesimal
CR-automorphisms} of $M$, denoted by $\frak{aut}_{CR}(M)$. The compact Lie group ${\sf
Aut}_{CR}(M)$, associated to this algebra is in fact the holomorphic symmetry group of $M$, the
local Lie
group of local biholomorphisms mapping $M$ to itself.
Determining such Lie algebras $\mathfrak{ aut}_{ CR} ( M)$
is the same as knowing the {\em
CR-symmetries} of $M$, a question which lies at the heart of the
(open) problem of classifying all local analytic CR-manifolds up to
biholomorphisms.

In our case of the universal CR-models, where it is assigned to each coordinate variable a so called {\it weight degree} ({\it cf.} \cite{Beloshapka2004,SHAM}), it can be plainly proved that (\cite{Beloshapka2004,SHAM}) the sought algebra $\frak{aut}_{CR}(M)$ of a CR-model $M$ is finite dimensional and takes
 the finite graded (in the sense of
Tanaka) form:
\begin{equation}
\label{g-int}
 \frak {aut}_{CR}(M)=\underbrace{\frak g_{-\rho}\oplus\cdots\oplus\frak
g_{-1}}_{\frak g_-}\oplus\frak g_0\oplus\underbrace{\frak
g_1\oplus\cdots\oplus\frak g_\varrho}_{\frak g_+} \ \ \ \ \
\rho,\,\varrho\in\mathbb N,
\end{equation}
where each component $\frak g_t$ is the Lie subalgebra of all
weighted homogeneous vector fields having the precise weight $t$. In this case the (negative) Lie subalgebra $\frak g_-$ is in fact the Levi-Tanaka algebra of $M$ and its associated compact Lie group $G_-$ is isomorphic to $M$, itself (\cite[Proposition 3]{Beloshapka2004}).

The classical approach of determining the already mentioned Lie algebras $\frak {aut}_{CR}(M)$ of each universal CR-model $M$ relies on constructing and solving some arising {\sc pde} system which has the role of the so-called {\it determining {\sc pde} system} in the general theory of Lie symmetries ({\it cf.} \cite{Beloshapka1997, Merker-Sabzevari-CEJM, 5-cubic, Mamai2009, Shananina2000}). The construction and solving these systems is quite complicated and time-consuming, in particular as much as the number of variables increases. Nevertheless, very recently in \cite{SHAM} we have designed a powerful algorithm enabling one to compute the desired algebras without constructing or solving any {\sc pde} system and just by employing some simple techniques from linear algebra and the modern theory of comprehensive Gr\"obner systems. This algorithm and its implementation in the {\sc Maple} software will be a great aid to proceed the computations pertinent to our second approach, will be discussed in section \ref{approach-2}.

\section{First approach:
Cartan equivalence problem}
\label{Cartan}

In this section, we attempt to employ the Cartan's algorithm for computing desired invariants of the moduli spaces. We explain this approach by computing the single invariant of the moduli space $\mathcal M(1,4)$. Our motivation behind employing this method in the case of the moduli spaces of universal CR-models came from $\S$12 of the long memoir \cite{5-cubic}, where it is studied the equivalence problem of totally nondegenerate CR-submanifolds in $\mathbb C^{4}$ to their unique universal CR-model, namely the 5-dimensional CR-cubic of the type $(1,3)$. Launching the classical algorithm of Cartan, first we need to construct an initial frame on the complexification of the under consideration CR-models.

\subsection{Constructing a frame}

For two integers $\bf a\in\mathbb C$ and $\bf b\in\mathbb R$, let $M({\bf a,b})\subset\mathbb C^{1+4}=\mathbb C\{z,w_1,w_2,w_3,w_4\}$, with $z=x+iy$ and $w_k=u_k+iv_k,\, k=1,\ldots,4$, be a six (real) dimensional Beloshapka's CR-model of CR-dimension 1 and codimension 4 represented as the graph of four defining polynomials:
\begin{equation}
\label{model-k=4}
\aligned
\Xi_1&:=w_1-\overline w_1-2i\,z\overline z=0,
\\
\Xi_2&:=w_2-\overline w_2-2i\,(z^2\overline z+z\overline z^2)=0,
\\
\Xi_3&:=w_3-\overline w_3-2\,(z^2\overline z-z\overline z^2)=0,
\\
\Xi_4&:=w_4-\overline w_4-2i\,({\bf a}\,z^3\overline z+\overline{\bf a}\,z\overline z^3)-2i{\bf b}\,z^2\overline z^2=0.
\endaligned
\end{equation}
Each pair $({\bf a,b})\in\mathbb C\times\mathbb R$ represents a distinct universal CR-model $M({\bf a,b})$ of the fixed type $(1,4)$. Without loss of generality, {\it one can assume throughout this paper that ${\bf a}$ is a nonzero complex integer}. Indeed, the case $({\bf a,b})=(0,0)$ is discarded since in this case, one easily checks that $M(0,0)$ is not totally nondegenerated. Moreover, it is easy to verify that all CR-models $M(0,{\bf b})$ with ${\bf b}\neq 0$ are holomorphically equivalent to the unique model $M(0,1)$ throughout the very simple transformation $w_4\mapsto {\bf b}w_4$. Let us denote by $\mathcal C(1,4)$ the class of these universal CR-models.

 According to \cite{Boggess-1991, Merker-Pocchiola-Sabzevari-5-CR-II, Merker-Porten-2006} and for a fixed model $M:=M({\bf a,b})$ of $\mathcal C(1,4)$, the $(1,0)$-complex tangent plane $T^{1,0}M$ is spanned by the single $(1,0)$-vector field:
\[
\mathcal L:=\frac{\partial}{\partial z}+A_1\,\frac{\partial}{\partial w_1}+A_2\,\frac{\partial}{\partial w_2}+A_3\,\frac{\partial}{\partial w_3}+A_4\,\frac{\partial}{\partial w_4},
\]
satisfying:
\[
\mathcal L |_{\Xi_k}\equiv 0, \,\,\,\,\,\,\, {\scriptstyle (k=1,\ldots,4)}.
\]
Then, applying this equality on the above four defining polynomials $\Xi_1,\ldots,\Xi_4$ and computing the coefficients $A_\bullet$ give:
\[
\mathcal L:=\frac{\partial}{\partial z}+2i\overline z\,\frac{\partial}{\partial w_1}+(4i\,z\overline z+2i\,\overline z^2)\,\frac{\partial}{\partial w_2}+(4\,z\overline z-2\,\overline z^2)\,\frac{\partial}{\partial w_3}+(6i{\bf a}\,z^2\overline z+2i\overline{\bf a}\,\overline z^3+4i{\bf b}\,z\overline z^2)\,\frac{\partial}{\partial w_4}.
\]
Here, the expression of $\mathcal L$ is presented as a
vector field which lives in a neighborhood of $M$ in
$\mathbb C^5$, while $M$ itself, is a real 6-dimensional
surface equipped with the six real coordinates
$x,y,u_1,u_2,u_3,u_4$. Thus, in order to express $\mathcal L$
{\it intrinsically}, one must drop $\frac{\partial}{\partial
v_k}$ for $k=1,\ldots,4$ and also simultaneously replace each $v_k$ by its expression in
\thetag{\ref{model-k=4}}.
Then, after expanding $\mathcal L$ in real and
imaginary parts one gets:
\[
\mathcal L:=\frac{\partial}{\partial z}+i\overline z\,\frac{\partial}{\partial u_1}+(2i\,z\overline z+i\,\overline z^2)\,\frac{\partial}{\partial u_2}+(2\,z\overline z-\overline z^2)\,\frac{\partial}{\partial u_3}+(3i\,{\bf a}\,z^2\overline z+i\,\overline{\bf a}\,\overline z^3+2i{\bf b}\,z\overline z^2)\,\frac{\partial}{\partial u_4}.
\]
While this vector field generates the $(1,0)$-complex bundle $T^{1,0}M$, its conjugation:
\[
\overline{\mathcal L}:=\frac{\partial}{\partial \overline z}-i z\,\frac{\partial}{\partial u_1}-(2i\,z\overline z+i\,z^2)\,\frac{\partial}{\partial u_2}+(2\,z\overline z- z^2)\,\frac{\partial}{\partial u_3}-(3i\,\overline{\bf a}\,z\overline z^2+i\,{\bf a}\, z^3+2i{\bf b}\,z^2\overline z)\,\frac{\partial}{\partial u_4}
\]
is, as well, the single generator of the $(0,1)$-complex tangent bundle $T^{0,1}M=\overline{T^{1,0}}M$.

So far, two vector fields $\mathcal L$ and $\overline{\mathcal L}$ are in fact the first elements of the sought frame for the complexified bundle $\mathbb C\otimes TM$. But still we need $6-2=4$ more independent complex vector fields to finish the construction of this frame. The totally nondegeneracy property of the under consideration universal CR-models enables one to compute the remaining necessary vector fields as the iterated Lie brackets of $\mathcal L$ and $\overline{\mathcal L}$ ({\it see} \cite{Beloshapka2004, Mamai} for the precise definition of totally nondegeneracy). At first, let us compute the length two Lie bracket $[\mathcal L,\overline{\mathcal L}]$. By plain computation, one finds it as an imaginary vector field. In order to get a real one, let us multiply it by $i$:
\[
\mathcal T:=i\,[\mathcal L,\mathcal{\overline L}]=2\,\frac{\partial}{\partial u_1}+4\,(z+\overline z)\,\frac{\partial}{\partial u_2}-4i\,(z-\overline z)\,\frac{\partial}{\partial u_2}+(6{\bf a}\,z^2+6\overline{\bf a}\,\overline z^2+8{\bf b}\,z\overline z)\,\frac{\partial}{\partial u_4}.
\]
Next, computing the length three iterated brackets gives:
\[
\aligned
\mathcal S&:=[\mathcal L,\mathcal T]=4\,\frac{\partial}{\partial u_2}-4i\,\frac{\partial}{\partial u_3}+(12 {\bf a}\,z+8{\bf b}\,\overline z)\,\frac{\partial}{\partial u_4},
\\
\overline{\mathcal S}&:=[\overline{\mathcal L},\mathcal T]=4\,\frac{\partial}{\partial u_2}+4i\,\frac{\partial}{\partial u_3}+(12 \overline{\bf a}\,\overline z+8{\bf b}\, z)\,\frac{\partial}{\partial u_4}.
\endaligned
\]
At the moment, we need just one more independent vector field. For this, computing the length four Lie bracket $[\mathcal L,\mathcal S]$ brings the expression:
\begin{equation}
\label{[l,lbar]}
[\mathcal L,\mathcal S]=12{\bf a}\,\frac{\partial}{\partial u_4}.
 \end{equation}
  Due to the fact that ${\bf a}\neq 0$, then one can consider the last sought vector field as the {\it real} one:
\begin{equation*}
\label{U}
\mathcal U:=\frac{1}{\bf a}\,[\mathcal L,\mathcal S]=12\frac{\partial}{\partial u_4}.
\end{equation*}

A glance on the expressions of the already computed six vector fields shows that they are linearly independent and thus we have;

\begin{Lemma}
The six vector fields $\mathcal L,\overline{\mathcal L},\mathcal T,\mathcal S,\overline{\mathcal S},\mathcal U$ construct a (complex) frame for the complexification $TM\otimes\mathbb C$ of the tangent bundle of $M$.
\end{Lemma}

The Lie commutators of the computed vector fields are displayed in the following table:
\medskip
\begin{center}
\footnotesize
\begin{tabular} [t] { c | c c c c c c }
& ${\mathcal L}$ & $\overline{\mathcal L}$ & ${\mathcal T}$ & ${\mathcal S}$ & $\overline{\mathcal S}$
& ${\mathcal U}$
\\
\hline
\medskip
${\mathcal L}$ & $0$ & $-i{\mathcal T}$ & ${\mathcal S}$ & ${\bf a}{\mathcal U}$ & $\frac{2\bf b}{3}{\mathcal U}$ & $0$
\\
$\overline{\mathcal L}$ & $*$ & $0$ & $\overline{\mathcal S}$ & $\frac{2\bf b}{3}{\mathcal U}$ & $\overline{\bf a}{\mathcal U}$ & $0$
\\
${\mathcal T}$ & $*$ & $*$ & $0$ & $0$ & $0$ & $0$
\\
${\mathcal S}$ & $*$ & $*$ & $*$ & $0$ & $0$ & $0$
\\
$\overline{\mathcal S}$ & $*$ & $*$ & $*$ & $*$ & $0$ & $0$
\\
${\mathcal U}$ & $*$ & $*$ & $*$ & $*$ & $*$ & $0$
\end{tabular}
\end{center}

\subsection{Passage to a dual coframe and its Darboux-Cartan structure}

Now let us introduce the coframe:
\begin{center}
$\{\mu_0, \sigma_0,\overline\sigma_0,\rho_0,\zeta_0,\overline\zeta_0\}$ which is dual to the constructed frame $\{\mathcal U, \mathcal S,\overline{\mathcal S},\mathcal T,\mathcal L,\overline{\mathcal L}\}$.
\end{center}
What we need at this step of launching the Cartan's algorithm of equivalence is to initially know the expressions of the five 2-forms $d\mu_0, d\sigma_0,d\overline\sigma_0,d\rho_0,d\zeta_0,d\overline\zeta_0$ in terms of the wedge products of the original 1-forms $\mu_0, \sigma_0,\overline\sigma_0,\rho_0,\zeta_0,\overline\zeta_0$. For this aim, we need the following well-known duality correspondence;

\begin{Lemma}
\label{lem-d}
Given a frame $\big\{ \mathcal{ L}_1, \dots, \mathcal{
L}_n\big\}$ on an open subset of $\mathbb R^n$ enjoying the Lie structure:
\[
\big[\mathcal{L}_{i_1},\,\mathcal{L}_{i_2}\big] =
\sum_{k=1}^n\,a_{i_1,i_2}^k\,\mathcal{L}_k \ \ \ \ \ \ \ \ \ \ \ \ \
{\scriptstyle{(1\,\leqslant\,i_1\,<\,i_2\,\leqslant\,n)}},
\]
where the $a_{ i_1, i_2}^k$ are certain functions on $\mathbb R^n$, the dual
coframe $\{ \omega^1, \dots, \omega^n \}$ satisfying by definition:
\[
\omega^k
\big(\mathcal{L}_i\big)
=
\delta_i^k
\]
enjoys a quite similar Darboux-Cartan structure, up to an overall minus sign:
\[
d\omega^k = - \sum_{1\leqslant i_1<i_2\leqslant n}\,
a_{i_1,i_2}^k\,\omega^{i_1}\wedge\omega^{i_2}
\ \ \ \ \ \ \ \ \ \ \ \ \ {\scriptstyle{(k\,=\,1\,\cdots\,n)}}.
\]
\end{Lemma}

Thanks to this Lemma, minding the overall minus sign, we can readily
find the expressions of the exterior derivatives of our five $1$-forms
that provide the associated Darboux-Cartan structure:

\begin{equation}
\label{d0}
\aligned
d\mu_0&={\bf a}\,\sigma_0\wedge\zeta_0+\textstyle{\frac{2\bf b}{3}}\,\overline\sigma_0\wedge\zeta_0+\frac{2\bf b}{3}\,\sigma_0\wedge\overline\zeta_0+\overline{\bf a}\,\overline\sigma_0\wedge\overline\zeta_0,
\\
d\sigma_0&=\rho_0\wedge\zeta_0, \ \ \ \ \ \ d\overline\sigma_0=\rho_0\wedge\overline\zeta_0,
\\
d\rho_0&=i\,\zeta_0\wedge\overline\zeta_0,
\\
d\zeta_0&=0, \ \ \ \ \ \ \ \ \ \ \ \ \ \ \ \ d\overline\zeta_0=0.
\endaligned
\end{equation}

\subsection{Ambiguity matrix}

After providing the initial frame $\{\mathcal L,\overline{\mathcal L},\mathcal T,\mathcal S,\overline{\mathcal S},\mathcal U\}$, now we can think about constructing the so-called ambiguity matrix of the under consideration equivalence problem, what encodes the equivalences of two arbitrary elements of the CR-class $\mathcal C(1,4)$.

Suppose that two given arbitrary CR-models $M:=M({\bf a,b})$ and $M':=M({\bf a',b'})$ belonging to this class are CR-equivalent through some local biholomorphism:
\[\aligned
h:M&\longrightarrow M'
\\
(z,w)& \ \mapsto (z'(z,w),w(z,w)).
\endaligned
\]
Then, the associated differential of
$h$, namely:
\[
h_\ast
\colon\ \ \
TM
\longrightarrow
T{M'}
\]
induces a push-forward complexified map, still
denoted with the same symbol:
\[
h_\ast
\colon\ \ \
\mathbb{C}\otimes
TM\longrightarrow
\mathbb{C}
\otimes
T{M'},
\]
which is naturally defined by ({\it see} \cite[Subsection
3.1]{Boggess-1991}):
\[
\aligned
h_\ast\big({\sf z}\otimes_\mathbb R \mathcal X\big)
:=
{\sf z}\otimes_\mathbb R
h_\ast(\mathcal X), \ \ \ \ \ \ \ {\footnotesize {\sf z}\in \mathbb
C, \ \ \ \mathcal X\in T_pM},\ \ \
p\in M.
\endaligned
\]

Now, assume that our new equivalent CR-model $M'$ is equipped by the complex frame $\{\mathcal L',\overline{\mathcal L'},\mathcal T',\mathcal S',\overline{\mathcal S'},\mathcal U'\}$ with the expressions similar to those of $M$ in new coordinates $(z',w'_1,w'_2,w'_3,w'_4)$. Since any biholomorphic equivalence $h$ satisfies  ({\it see} \cite{Merker-Pocchiola-Sabzevari-5-CR-II}):
\[
h_\ast(T_p^{1,0}M)=T_{h(p)}^{1,0}M' \ \ \  \textrm{and} \ \ \ h_\ast(\overline{\mathcal Y})=\overline{h_\ast(\mathcal Y)}, \ \ \ \mathcal Y\in\mathbb C\otimes TM,
\]
then one accordingly deduces that $h$ maps the single generators of $T^{1,0}M$ and $T^{0,1}M$ as:
\begin{equation}
\label{h-L}
h_\ast(\mathcal L)=a_1\,\mathcal L' \ \ \ \textrm{and} \ \ \ h_\ast(\overline{\mathcal L})=\overline a_1\,\overline{\mathcal L'}
\end{equation}
for some {\it nonzero} complex functions $a_1$ defined on $M'$.

Next, let us look at what happens with Lie brackets. Since the
differential operator commutes with brackets, we have:
\[
\aligned h_{\ast}(\mathcal T)
= h_{\ast}\big(i[\mathcal
L,\overline{\mathcal L}]\big)=
i\,h_\ast\big([\mathcal L,\overline{\mathcal
L}]\big)
=i\,
\big[h_\ast(\mathcal L),h_\ast(\overline{\mathcal L})\big]
&
=i\big[a_1\mathcal L',\overline a_1\overline{\mathcal L}'\big],
\endaligned
\]
and by expanding this last bracket one obtains:
\begin{equation}
\label{h-T}
\aligned
h_{\ast}(\mathcal T)
&
=
a_1\overline a_1\cdot i\,
\big[\mathcal{L}',\,\overline{\mathcal{L}'}\big]
\underbrace{-
i\,\overline a_1\overline{\mathcal{L}}'( a_1)}_{
=:\,a_2}
\cdot
\mathcal{L}'
+
i\, a_1\,\mathcal{L}'\big(\overline a_1\big)
\cdot
\overline{\mathcal{L}}'
\\
&=: a_1\overline a_1\mathcal T'+ a_2\mathcal
L'+\overline a_2\overline{\mathcal L}',
\endaligned
\end{equation}
for some new appearing complex function $a_2$ defined on $M'$. Proceeding along the same lines of computations, one also finds the image of the next basis fields $\mathcal S, \overline{\mathcal U}$ and $\mathcal U$ as:
\begin{equation}
\label{h-S-U}
\aligned
h_\ast(\mathcal S)&:=a_1^2\overline a_1\,\mathcal S'+a_3\,\mathcal T'+a_4\mathcal L'+\overline a_5\,\overline{\mathcal L'},
\\
h_\ast(\overline{\mathcal S})&:=a_1\overline a_1^2\,\overline{\mathcal S'}+\overline a_3\,\mathcal T'+ a_5\mathcal L'+\overline a_4\,\overline{\mathcal L'},
\\
h_\ast(\mathcal U)&:=a_1^3\overline a_1\,\mathcal U'+a_6\,\mathcal S'+a_7\,\mathcal T'+a_8\mathcal L'+a_9\,\overline{\mathcal L'},
\endaligned
\end{equation}
for some complex functions $a_3,\ldots,a_9$ defined on $M'$.

\subsubsection{Still some more simplifications}

Having $\mathcal U$ as a real vector field implies on the one hand the equality $h_\ast(\overline{\mathcal U})=h_\ast(\mathcal U)$. On the other hand, due to the intrinsic properties of the diffeomorphism $h$ one has:
\[
\aligned
h_\ast(\overline{\mathcal U})=\overline{h_\ast(\mathcal U)}=a_1\overline a_1^3\,{\mathcal U'}+\overline a_6\,\overline{\mathcal S'}+\overline a_7\,\mathcal T'+\overline a_9\mathcal L'+\overline a_8\,\overline{\mathcal L'}.
\endaligned
\]
Now, it follows immediately after comparison the above equal expressions of $h_\ast(\mathcal U)$ and $h_\ast(\overline{\mathcal U})$ that:
\[
a_6\equiv 0, \,\,\,\, a_7=\overline a_7, \,\,\,\, a_9=\overline a_8.
\]
Thus, with minor modification on the indices one finds the following expression for $h_\ast(\mathcal U)$:
\begin{equation}
\label{h-U}
h_\ast(\mathcal U)=a_1^3\overline a_1\,\mathcal U'+a_6\,\mathcal T'+a_7\mathcal L'+\overline a_7\,\overline{\mathcal L'},
\end{equation}
for the functions $a_1$ and $a_7$ as above and for some {\it real-valued} function $a_6$. Still there is
another useful fact comes from comparing the coefficients of the real vector field $\mathcal U'$ in the both sides of the equality $\overline{h_\ast(\mathcal U)}=h_\ast(\mathcal U)$. Accordingly we have $a_1\overline a_1^3=a_1^3\overline a_1$ which immediately implies that  $a_1^3\overline a_1$ is {\it real} and

\begin{Lemma}
\label{a1}
The nonzero complex function $a_1$ enjoys the equality:
 \begin{equation*}
 a_1^2=\overline a_1^2
\end{equation*}
which equivalently means that $a_1$ is either real or imaginary,
\end{Lemma}
\noindent
a fact that will have a great influence to simplify the next computations throughout this section.

Summing up the results and according to the expressions \thetag{\ref{h-L}} -- \thetag{\ref{h-U}}, one finds out that there exists a local biholomorphism $h$ between two arbitrary CR-models $M,M'\in\mathcal C(1,4)$ {\it if and only if} the associated moving frames can be mapped to each other as follows through some complex functions $a_1,\ldots,a_7$:
\[\footnotesize
\left(
  \begin{array}{c}
    \mathcal U \\
    \mathcal S \\
    \overline{\mathcal S} \\
    \mathcal T \\
    \mathcal L \\
    \overline{\mathcal L} \\
  \end{array}
\right)
=
\left(
  \begin{array}{cccccc}
    a_1^3\overline a_1 & 0 & 0 & a_6 & a_7 & \overline a_7 \\
    0 & a_1^2\overline a_1 & 0 & a_3 & a_4 & \overline a_5 \\
    0 & 0 & a_1\overline a_1^2 & \overline a_3 & a_5 & \overline a_4 \\
    0 & 0 & 0 & a_1\overline a_1 & a_2 & \overline a_2 \\
    0 & 0 & 0 & 0 & a_1 & 0 \\
    0 & 0 & 0 & 0 & 0 & \overline a_1 \\
  \end{array}
\right).
\left(
  \begin{array}{c}
    \mathcal U' \\
    \mathcal S' \\
    \overline{\mathcal S'} \\
    \mathcal T' \\
    \mathcal L' \\
    \overline{\mathcal L'} \\
  \end{array}
\right).
\]

\subsection{Setting up the equivalence problem}

According to the general principles (\cite{Olver-1995}),
the so-called {\sl lifted coframe} in terms
of the dual basis of $1$-form then becomes,
after a plain matrix transposition:

\begin{equation}
\label{lifted-1}
\footnotesize
\left(
  \begin{array}{c}
    \mu \\
    \sigma \\
    \overline{\sigma} \\
    \rho \\
    \zeta \\
    \overline{\zeta} \\
  \end{array}
\right)
=
\underbrace{\left(
  \begin{array}{cccccc}
    a_1^3\overline a_1 & 0 & 0 & 0 & 0 & 0 \\
    0 & a_1^2\overline a_1 & 0 & 0 & 0 & 0 \\
    0 & 0 & a_1\overline a_1^2 & 0 & 0 & 0 \\
    a_6 & a_3 & \overline a_3 & a_1\overline a_1 & 0 & 0 \\
    a_7 & a_4 & a_5 & a_2 & a_1 & 0 \\
    \overline a_7 & \overline a_5 & \overline a_4 & \overline a_2 & 0 & \overline a_1 \\
  \end{array}
\right)}_{g}.
\left(
  \begin{array}{c}
    \mu_0 \\
    \sigma_0 \\
    \overline{\sigma_0} \\
    \rho_0 \\
    \zeta_0 \\
    \overline{\zeta_0} \\
  \end{array}
\right).
\end{equation}
In this case, the collection:
\[
G:=\big\{g: \,\, a_6\in\mathbb R, \ \ a_k\in\mathbb C\,\, \textrm{for}\, k\neq 6 \,\, \textrm{and}\,\, a_1^2=\overline a_1^2\big\}\subset\sf{GL}(6)
\]
 constitutes a matrix Lie group which is called the {\sl structure group} of the equivalence problem.

\subsection{Structure equations}

Let $\Theta:=(\mu,\sigma,\overline\sigma,\rho,\zeta,\overline\zeta)^t$ and $\Theta_0:=(\mu_0,\sigma_0,\overline\sigma_0,\rho_0,\zeta_0,\overline\zeta_0)^t$ be the lifted and initial coframes as above. Then the lifted coframe structure \thetag{\ref{lifted-1}} can be read as:
\[
\Theta=g\cdot\Theta_0.
\]
Differentiating this fundamental equality gives:
\begin{equation}
\label{struc-eq-preparation}
{\sf d}\Theta={\sf d}g\wedge\Theta_0+g\cdot{\sf d}\Theta_0.
\end{equation}
The first term in this expression of ${\sf d}\Theta$ can be written in the form:
\[
{\sf d}g\wedge\Theta_0=\underbrace{{\sf d}g\cdot g^{-1}}_{=:\omega_{\sf MC}}\wedge \underbrace{g\cdot\Theta_0}_{\Theta},
\]
in which after necessary computations and summing with  $g\cdot {\sf d}\Theta_0$ one receives:
\begin{equation}
\footnotesize\aligned
\label{MC-4-1}
{\sf d}\Theta={\sf d}\left(
                   \begin{array}{c}
                     \mu \\
                     \sigma \\
                     \overline\sigma \\
                     \rho \\
                     \zeta \\
                     \overline\zeta \\
                   \end{array}
                 \right)&=\left(
                   \begin{array}{cccccc}
                     3\,\alpha_1+\overline\alpha_1 & 0 & 0 & 0 & 0 & 0 \\
                     0 & 2\,\alpha_1+\overline\alpha_1 & 0 & 0 & 0 & 0 \\
                     0 & 0 & \alpha_1+2\,\overline\alpha_1 & 0 & 0 & 0 \\
                     \alpha_2 & \alpha_3 & \overline\alpha_3 & \alpha_1+\overline\alpha_1 & 0 & 0 \\
                     \alpha_4 & \alpha_5 & \alpha_6 & \alpha_7 & \alpha_1 & 0 \\
                     \overline\alpha_4 & \overline\alpha_6 & \overline\alpha_5 & \overline\alpha_7 & 0 & \overline\alpha_1 \\
                   \end{array}
                 \right)
\wedge
                 \left(
                   \begin{array}{c}
                     \mu \\
                     \sigma \\
                     \overline\sigma \\
                     \rho \\
                     \zeta \\
                     \overline\zeta \\
                   \end{array}
                 \right)
\\
&\ \ \ \ \ \ \ \ \ \ \ \ +
                              \left(
  \begin{array}{c}
    a_1^3\overline a_1\,{\sf d}\mu_0 \\
     a_1^2\overline a_1\,{\sf d}\sigma_0 \\
     a_1\overline a_1^2\,{\sf d}\overline\sigma_0 \\
    a_6\,{\sf d}\mu_0+a_3\,{\sf d}\sigma_0+\overline a_3\,{\sf d}\overline\sigma_0+a_1\overline a_1 {\sf d}\rho_0 \\
    a_7\,{\sf d}\mu_0+a_4\,{\sf d}\sigma_0+a_5\,{\sf d}\overline\sigma_0+a_2\,{\sf d}\rho_0+a_1\,{\sf d}\zeta_0 \\
    \overline a_7\,{\sf d}\mu_0+\overline a_5\,{\sf d}\sigma_0+\overline a_4\,{\sf d}\overline\sigma_0+\overline a_2\, {\sf d}\rho_0+\overline a_1\,{\sf d}\overline\zeta_0 \\
  \end{array}
\right).
\endaligned
\end{equation}
 Here, $\omega_{\sf MC}$ is the {\sl Maurer-Cartan matrix} associated to our structure Lie group $G$ constituted by the following Maurer-Cartan 1-forms\footnote{These expressions are a bit simplified according to Lemma \ref{a1}} $\alpha_\bullet$:
\[\aligned
\alpha_1&:=\textstyle\frac{1}{a_1}\,{\sf d}a_1,
\\
\alpha_2&:=\textstyle\frac{1}{a_1^3\overline a_1}\,{\sf d}a_6-\frac{a_6}{a_1^4\overline a_1}\,{\sf d} a_1-\frac{a_6}{a_1\overline a_1^4}\,{\sf d}\overline a_1,
\\
\alpha_3&:=\textstyle\frac{1}{a_1^2\overline a_1}\,{\sf d}a_3-\frac{a_3}{a_1^3\overline a_1}\,{\sf d}a_1-\frac{a_3}{a_1^2\overline a_1^2}\,{\sf d}\overline a_1,
\\
\alpha_4&:=\textstyle\frac{1}{a_1^3\overline a_1}\,{\sf d}a_7-\frac{a_6}{a_1^4\overline a_1^2}\,{\sf d}a_2-\frac{a_1a_7\overline a_1-a_2a_6}{a_1^5\overline a_1^2}\,{\sf d}a_1,
\\
\alpha_5&:=\textstyle\frac{1}{a_1^2\overline a_1}\,{\sf d}a_4-\frac{a_3}{a_1^3\overline a_1^2}\,{\sf d}a_2-\frac{a_1a_4\overline a_1-a_2a_3}{a_1^4\overline a_1^2}\,{\sf d}a_1,
\\
\alpha_6&:=\textstyle\frac{1}{a_1\overline a_1^2}\,{\sf d}a_5-\frac{\overline a_3}{a_1^2\overline a_1^3}\,{\sf d}a_2-\frac{a_1a_5\overline a_1-a_2\overline a_3}{a_1^3\overline a_1^3}\,{\sf d}a_1,
\\
\alpha_7&:=\textstyle \frac{1}{a_1\overline a_1}\,{\sf d}a_2-\frac{a_2}{a_1^2\overline a_1}\,{\sf d}a_1.
\endaligned
\]
Thanks to  Lemma \ref{a1}, here $\alpha_2$ is a {\it real} 1-form and:
\[
\alpha_1=\frac{{\sf d}a_1}{a_1}=\frac{{\sf d}\overline a_1}{\overline a_1}=\overline\alpha_1,
\]
which means that $\alpha_1$ is {\it real}, too. Then, one can replace $\overline\alpha_1$ by $\alpha_1$ in the structure equations \thetag{\ref{MC-4-1}}.

In the expression of ${\sf d}\Theta$ in \thetag{\ref{MC-4-1}}, the second terms, namely the last $6\times 1$ matrix, is the only part which still includes the initial 2-forms of ${\sf d}\Theta_0$. However, it is possible to express the 2-forms ${\sf d}\mu_0,{\sf d}\sigma_0,\ldots,{\sf d}\overline\zeta_0$ as some certain combinations of the wedge products of $\mu,\sigma,\overline\sigma,\rho,\zeta,\overline\zeta$ by regarding first the equality $\Theta_0=g^{-1}\cdot\Theta$, namely:
\begin{equation}
\label{mu0--mu}
\aligned
\mu_0&=\textstyle\frac{1}{a_1^3\overline a_1}\,\mu,
\\
\sigma_0&=\textstyle\frac{1}{a_1^2\overline a_1}\,\sigma,
\\
\rho_0&=-\textstyle\frac{a_6}{a_1^4\overline a_1^2}\,\mu-\frac{a_3}{a_1^3\overline a_1^2}\sigma-\frac{\overline a_3}{a_1^2\overline a_1^3}\,\overline\sigma+\frac{1}{a_1\overline a_1}\,\rho,
\\
\zeta_0&=-\textstyle\frac{a_1a_7\overline a_1-a_2a_6}{a_1^5\overline a_1^2}\,\mu-\frac{a_1a_4\overline a_1-a_2a_3}{a_1^4\overline a_1^2}\,\sigma-\frac{a_1a_5\overline a_1-a_2\overline a_3}{a_1^3\overline a_1^3}\,\overline\sigma-\frac{a_2}{a_1^2\overline a_1}\,\rho+\frac{1}{a_1}\,\zeta
\endaligned
\end{equation}
and next substituting it into the expressions of ${\sf d}\mu_0,{\sf d}\sigma_0,{\sf d}\overline\sigma_0,{\sf d}\rho_0,{\sf d}\zeta_0,{\sf d}\overline\zeta_0$ in \thetag{\ref{d0}}. Inserting the obtained expressions in \thetag{\ref{MC-4-1}} and taking it into account the equality $\alpha_1=\overline\alpha_1$, one finally receives the so-called {\sl structure equations} of the problem as:
\begin{equation}
\footnotesize\aligned
\label{struc-eq-1}
{\sf d}\mu&:=4\,\alpha_1\wedge\mu+
 \\
&+T_1\,\mu\wedge\sigma+\overline T_1\,\mu\wedge\overline\sigma+ T_2\,\sigma\wedge\overline\sigma+T_3\,\sigma\wedge\rho+{\bf a}\,\sigma\wedge\zeta+{\sf c}\,\sigma\wedge\overline\zeta+
\\
& \ \ \ \ \ \ \ \ \ \ \ \ \ \ \ \ \ \ \ \ \ \ \ \ \ \ \ \ \ \ \ \ \ \ \ \ \ \ \ \ \ \ \ \ \ \ \ \ \ \ \ \ \ \ \ \ \ +\overline T_3\,\overline\sigma\wedge\rho+ {\sf c}\,\overline\sigma\wedge\zeta+\overline{\bf a}\,\overline\sigma\wedge\overline\zeta,
\\
{\sf d}\sigma&:=3\,\alpha_1\wedge\sigma+
\\
&+ U_1\,\mu\wedge\sigma+U_2\,\mu\wedge\overline\sigma+U_3\,\mu\wedge\rho+U_4\,\mu\wedge\zeta+
\\
&\ \ \ \ \ \ \ \ \ \ \ \ \ \ \ \ \ \ \  +U_5\,\sigma\wedge\overline\sigma+U_6\,\sigma\wedge\rho+U_7\,\sigma\wedge\zeta+
\\
&
\ \ \ \ \ \ \ \ \ \ \ \ \ \ \ \ \  \ \ \ \ \ \ \ \ \ \ \ \ \ \ \ \ \ \ \ \ \ +U_8\,\overline\sigma\wedge\rho+\overline U_7\,\overline\sigma\wedge\zeta+\rho\wedge\zeta,
\\
{\sf d}\rho&:=\alpha_2\wedge\mu+\alpha_3\wedge\sigma+\overline\alpha_3\wedge\overline\sigma+2\,\alpha_1\wedge\rho+
\\
& +V_1\,\mu\wedge\sigma+\overline V_1\,\mu\wedge\overline\sigma+V_2\,\mu\wedge\rho+V_3\,\mu\wedge\zeta+\overline V_3\,\mu\wedge\overline\zeta+
\\
&\ \ \ \ \ \ \ \ \ \ \ \ \ \ \ \ \ \ \ +V_4\,\sigma\wedge\overline\sigma+V_5\,\sigma\wedge\rho+V_6\,\sigma\wedge\zeta+V_7\,\sigma\wedge\overline\zeta+
\\
& \ \ \ \ \ \ \ \ \ \ \ \ \ \ \ \  \ \ \ \ \ \ \ \ \ \ \ \ \ \ \ \  \ \ \ \ \ +\overline V_5\,\overline\sigma\wedge\rho+\overline V_7\,\overline\sigma\wedge\zeta+\overline V_6\,\overline\sigma\wedge\overline\zeta+
\\
& \ \ \ \ \ \ \ \ \ \ \ \ \ \ \ \  \ \ \ \ \ \ \ \ \ \ \ \ \ \ \ \  \ \ \ \ \ \ \ \ \ \ \ \ \ \ \ \ \ \ \ \ \ \ \ \ \ \ +V_8\,\rho\wedge\zeta+\overline V_8\,\rho\wedge\overline\zeta+i\,\zeta\wedge\overline\zeta,
\\
{\sf d}\zeta&:=\alpha_4\wedge\mu+\alpha_5\wedge\sigma+\alpha_6\wedge\overline\sigma+\alpha_7\wedge\rho+\alpha_1\wedge\zeta+
\\
&+W_1\,\mu\wedge\sigma+W_2\,\mu\wedge\overline\sigma+W_3\,\mu\wedge\rho+W_4\,\mu\wedge\zeta+W_5\,\mu\wedge\overline\zeta+
\\
& \ \ \ \ \ \ \ \ \ \ \ \ \ \ \ \ \ \ \ \ +W_6\,\sigma\wedge\overline\sigma+W_7\,\sigma\wedge\rho+W_8\,\sigma\wedge\zeta+W_9\,\sigma\wedge\overline\zeta+
\\
&\ \ \ \ \ \ \ \ \ \ \ \ \ \ \ \ \ \ \ \  \ \ \ \ \ \ \ \ \ \ \ \ \ \ \ \ \ \ \ +W_{10}\,\overline\sigma\wedge\rho+W_{11}\,\overline\sigma\wedge\zeta+W_{12}\,\overline\sigma\wedge\overline\zeta+
\\
&\ \ \ \ \ \ \ \ \ \ \ \ \ \ \ \ \ \ \ \  \ \ \ \ \ \ \ \ \ \ \ \ \ \ \ \ \ \ \ \ \ \ \ \ \ \ \ \ \ \ \ \ \ \ +W_{13}\,\rho\wedge\zeta+W_{14}\,\rho\wedge\overline\zeta+W_{15}\,\zeta\wedge\overline\zeta,
\endaligned
\end{equation}
fully expressed in terms of the wedge products of the Maurer-Cartan and lifted 1-forms and without remaining any initial form among them.
 Saving the space, here we do not present the explicit (some of them long) expressions of the appearing {\sl torsion coefficients} $T_\bullet,U_\bullet,V_\bullet,W_\bullet$, though one can find them in the {\sc Maple} worksheet \cite{Maple}. Nevertheless, we present each of these expressions as soon as we require it in the next steps. One just notices that a glance on these expressions reveals that amongst them, $V_2$ is a {\it real function} while $T_2$ and $V_4$ are {\it imaginary}. Moreover, here the coefficient $\sf c$ of $\sigma\wedge\overline\zeta$ in the expression of ${\sf d}\mu$ is equal to $\frac{2}{3}.\frac{a_1\bf b}{\overline a_1}$, where according to Lemma \ref{a1} the quotient $\frac{a_1}{\overline a_1}$ is either $1$ or $-1$. Hence we have:
\begin{equation}
\label{c}
{\sf c}=\textstyle\frac{2}{3}\,{\bf b} \ \ \ \ \textrm{or} \ \ \ \ {\sf c}=-\frac{2}{3}\,{\bf b}.
\end{equation}

\subsection{First-loop absorbtion and normalization}

After providing the preliminary requirements, now we are ready to start the process of Cartan's algorithm. As the first step, we have to apply the absorbtion and normalization parts based on the results of Proposition \ref{prop-changes}. Accordingly, let us replace each Maurer-Cartan form $\alpha_\bullet$ by:
\[
\alpha_k\mapsto \alpha_k+p_k\, \mu+q_k\,\sigma+r_k\, \overline\sigma+s_k\,\rho+t_k\, \zeta+u_k\, \overline\zeta, \ \ \ \ \ \ \ \ \  {\scriptstyle (k=1,\ldots,7)},
\]
for some arbitrary coefficient-functions $p_k,q_k,r_k,s_k,t_k,u_k$.
These substitutions convert the structure equations \thetag{\ref{struc-eq-1}} into the form:
\begin{equation*}
\footnotesize\aligned
{\sf d}\mu&:=4\,\alpha_1\wedge\mu+
 \\
&+(T_1-4\,q_1)\,\mu\wedge\sigma+(\overline T_1-4\,r_1)\,\mu\wedge\overline\sigma-4s_1\,\mu\wedge\rho-4t_1\,\mu\wedge\zeta-4u_1\,\mu\wedge\overline \zeta+
\\
&
+ T_2\,\sigma\wedge\overline\sigma+T_3\,\sigma\wedge\rho+{\bf a}\,\sigma\wedge\zeta+{\sf c}\,\sigma\wedge\overline\zeta +\overline T_3\,\overline\sigma\wedge\rho+ {\sf c}\,\overline\sigma\wedge\zeta+\overline{\bf a}\,\overline\sigma\wedge\overline\zeta,
\\
{\sf d}\sigma&:=3\,\alpha_1\wedge\sigma+
\\
&+ (U_1+3\,p_1)\,\mu\wedge\sigma+U_2\,\mu\wedge\overline\sigma+U_3\,\mu\wedge\rho+U_4\,\mu\wedge\zeta+
\\
&  +(U_5-3\,r_1)\,\sigma\wedge\overline\sigma+(U_6-3\,s_1)\,\sigma\wedge\rho+(U_7-3\,t_1)\,\sigma\wedge\zeta -3u_1\,\sigma\wedge\overline\zeta+U_8\,\overline\sigma\wedge\rho+\overline U_7\,\overline\sigma\wedge\zeta+\rho\wedge\zeta,
\\
{\sf d}\rho&:=\alpha_2\wedge\mu+\alpha_3\wedge\sigma+\overline\alpha_3\wedge\overline\sigma+2\,\alpha_1\wedge\rho+
\\
& +(V_1-q_2+p_3)\,\mu\wedge\sigma+(\overline V_1+\overline p_3-r_2)\,\mu\wedge\overline\sigma+(V_2-s_2+2p_1)\,\mu\wedge\rho+(V_3-t_2)\,\mu\wedge\zeta+(\overline V_3-u_2)\,\mu\wedge\overline\zeta+
\\
&\ \ \  +(V_4-r_3+\overline r_3)\,\sigma\wedge\overline\sigma+(V_5-s_3+2q_1)\,\sigma\wedge\rho+(V_6-t_3)\,\sigma\wedge\zeta+(V_7-u_3)\,\sigma\wedge\overline\zeta+
\\
& +(\overline V_5-\overline s_3+2r_1)\,\overline\sigma\wedge\rho+(\overline V_7-\overline u_3)\,\overline\sigma\wedge\zeta+(\overline V_6-\overline t_3)\,\overline\sigma\wedge\overline\zeta +(V_8-2t_1)\,\rho\wedge\zeta+(\overline V_8-2u_1)\,\rho\wedge\overline\zeta+i\,\zeta\wedge\overline\zeta,
\\
{\sf d}\zeta&:=\alpha_4\wedge\mu+\alpha_5\wedge\sigma+\alpha_6\wedge\overline\sigma+\alpha_7\wedge\rho+\alpha_1\wedge\zeta+
\\
&+(W_1-q_4+p_5)\,\mu\wedge\sigma+(W_2-r_4+p_6)\,\mu\wedge\overline\sigma+(W_3-s_4+p_7)\,\mu\wedge\rho+(W_4-t_4+p_1)\,\mu\wedge\zeta+
(W_5-u_4)\,\mu\wedge\overline\zeta+
\\
& \ \ \ \ \ \ \ \ \ \ \ \ \ \ \ \ \ \ \ \ +(W_6-r_5+q_6)\,\sigma\wedge\overline\sigma+(W_7-s_5+q_7)\,\sigma\wedge\rho+(W_8-t_5+q_1)\,\sigma\wedge\zeta+
(W_9-u_5)\,\sigma\wedge\overline\zeta+
\\
&\ \ \ \ \ \ \ \ \ \ \ \ \ \ \ \ \ \ \ \  \ \ \ \ \ \ \ \ \ \ \ \ \ \ \ \ \ \ \ +(W_{10}-s_6+r_7)\,\overline\sigma\wedge\rho+(W_{11}-t_6+r_1)\,\overline\sigma\wedge\zeta+(W_{12}-u_6)\,\overline\sigma\wedge\overline\zeta+
\\
&\ \ \ \ \ \ \ \ \ \ \ \ \ \ \ \ \ \ \ \  \ \ \ \ \ \ \ \ \ \ \ \ \ \ \ \ \ \ \ \ \ \ \ \ \ \ \ \ \ \ \ \ \ \ +(W_{13}-t_7+s_1)\,\rho\wedge\zeta+(W_{14}-u_7)\,\rho\wedge\overline\zeta+(W_{15}-u_1)\,\zeta\wedge\overline\zeta.
\endaligned
\end{equation*}
In the situation that all $p_k,q_k,r_k,s_k,t_k,u_k$ are regarded as {\it arbitrary} coefficient-functions, one can determine them in such  a way that the appearing coefficients either vanish or being  equal to a constant integer as much as possible, the procedure which is called by {\sl normalization}. For instance, with $s_1=u_1=t_1\equiv 0$ and $q_1=\textstyle\frac{1}{4}T_1, r_1=\frac{1}{4}\overline T_1$, the coefficients of $\mu\wedge\sigma, \mu\wedge\overline\sigma,\mu\wedge\rho, \mu\wedge\zeta$ and $\mu\wedge\overline\zeta$ in the first expression ${\sf d}\mu$ vanish, identically. Nevertheless, to normalize as much as possible the appearing coefficients in the above expressions and besides appropriate determination of the coefficient-functions $p_k,q_k,r_k,s_k,t_k,u_k$, one also has to normalize torsion coefficients:
\begin{equation}
\label{essen-torsions}
\footnotesize\aligned
T_2&=\textstyle-\frac{1}{3\,a_1^2\overline a_1^4}\big(-3\overline{\bf a}\,a_1^2\overline a_1\overline a_5+3\,\overline{\bf a}\,a_1a_3\overline a_2+2{\bf b}\,a_1^2\overline a_1\overline a_4-2{\bf b}\,a_1\overline a_2\overline a_3-2{\bf b}\,a_1a_4\overline a_1^2+2{\bf b}\,a_2a_3\overline a_1+
\\
&+3{\bf a}\,a_1a_5\overline a_1^2-3{\bf a}\,a_2\overline a_1\overline a_3\big),
\\
T_3&=\textstyle-\frac{1}{3\,a_1\overline a_1^2}\,\big(2{\bf b}\,a_1\overline a_2+3{\bf a}\,a_2\overline a_1\big),
\\
U_2&=\frac{1}{a_1^4\overline a_1^3}\,\big(a_5a_6-a_7\overline a_3\big),
\\
U_3&=\frac{1}{a_1^3\overline a_1}\,a_7, \ \ \
U_4=-\frac{1}{a_1^3\overline a_1}\,a_6, \ \ \
U_6=\frac{1}{a_1^2\overline a_1}\,a_4, \ \ \
U_7=-\frac{1}{a_1^2\overline a_1}\,a_3, \ \ \
U_8=\frac{1}{a_1\overline a_1^2}\,a_5,
\\
U_5-\textstyle\frac{3}{4}\,\overline T_1&=-\frac{1}{4\,a_1^3\overline a_1^6}\big(4\,a_4 \overline a_1^3\overline a_3-4\,a_3a_5\overline a_1^3+3\overline{\bf a}\,a_1^4\overline a_1\overline a_7-3\overline{\bf a}\,a_1^3a_6 \overline a_2+2{\bf b}\,a_1a_7\overline a_1^4-2{\bf b}\,a_2a_6\overline a_1^3\big),
\\
V_8&=\frac{1}{a_1^2\overline a_1}\,\big(i\,a_1\overline a_2+a_3\big),
\\
W_{15}&=\frac{i}{a_1\overline a_1}\,a_2.
\endaligned
\end{equation}
A careful look at the expressions of these essential torsion coefficients shows that they can be normalized to zero when one puts:
\[
a_2=a_3=a_4=a_5=a_6=a_7\equiv 0.
\]
By such determination of the group parameters, the only remaining group parameter is $a_1$ and hence the structure equations \thetag{\ref{struc-eq-1}} take the simple form:
\begin{equation}
\aligned
\label{struc-eq-2}
{\sf d}\mu&:=4\,\alpha\wedge\mu+
{\bf a}\,\sigma\wedge\zeta+{\sf c}\,\sigma\wedge\overline\zeta+
{\sf c}\,\overline\sigma\wedge\zeta+\overline{\bf a}\,\overline\sigma\wedge\overline\zeta,
\\
{\sf d}\sigma&:=3\,\alpha\wedge\sigma+\rho\wedge\zeta,
\\
{\sf d}\rho&:=2\,\alpha\wedge\rho+i\,\zeta\wedge\overline\zeta,
\\
{\sf d}\zeta&:=\alpha\wedge\zeta,
\endaligned
\end{equation}
for $\alpha:=\frac{{\sf d}a_1}{a_1}$ and $\sf c$ as \thetag{\ref{c}}.

Applying any other normalization procedure is useless and thus one has to start the prolongation step. Moreover, as a consequence of the Cartan's Lemma ({\it see} \cite[Exercise 1.33]{Olver-1995}) the single remaining Maurer-Cartan form $\alpha$ is the unique 1-form satisfying the above structure equations.

\subsection{Prolongation}

Once the six group parameters $a_2,\ldots,a_7$ normalized to zero, the original structure group $G$ is reduced to the 1-dimensional matrix Lie group $G^{\sf red}$ of the matrices:
\[
\footnotesize\aligned
\left(
  \begin{array}{cccccc}
    a_1^3\overline a_1 & 0 & 0 & 0 & 0 & 0 \\
    0 & a_1^2\overline a_1 & 0 & 0 & 0 & 0 \\
    0 & 0 & a_1\overline a_1^2 & 0 & 0 & 0 \\
    0 & 0 & 0 & a_1\overline a_1 & 0 & 0 \\
    0 & 0 & 0 & 0 & a_1 & 0 \\
    0 & 0 & 0 & 0 & 0 & \overline a_1 \\
  \end{array}
\right),
\endaligned
\]
where, as before, $a_1$ is either a real or an imaginary parameter. According to principles of Cartan theory ({\it see} \cite[Proposition 12.1]{Olver-1995}), the equivalence problem of 6-dimensional CR-models $M({\bf a,b})$ can now be characterized in terms of the equivalence problem of the so-called {\sl prolonged spaces} $M({\bf a,b})\times G^{\sf red}$, equipped with the lifted coframe $(\mu,\sigma,\overline\sigma,\rho,\zeta,\overline\zeta,\alpha)$. One finds easily the structure equations associated to this new problem just by adding the exterior differentiation of the\,\,---\,\,formerly Maurer-Cartan and now initial\,\,---\,\,1-form $\alpha$, namely ${\sf d}\alpha={\sf d}\big(\frac{{\sf d}a_1}{a_1}\big)=0$, to the structure equations \thetag{\ref{struc-eq-2}}. Hence, the final structure equations of this new prolonged space readily manifest themselves as:

\begin{equation}
\aligned
\label{struc-eq-prolonged}
{\sf d}\mu&:=4\,\alpha\wedge\mu+
{\bf a}\,\sigma\wedge\zeta+{\sf c}\,\sigma\wedge\overline\zeta+
{\sf c}\,\overline\sigma\wedge\zeta+\overline{\bf a}\,\overline\sigma\wedge\overline\zeta,
\\
{\sf d}\sigma&:=3\,\alpha\wedge\sigma+\rho\wedge\zeta,
\\
{\sf d}\rho&:=2\,\alpha\wedge\rho+i\,\zeta\wedge\overline\zeta,
\\
{\sf d}\zeta&:=\alpha\wedge\zeta,
\\
{\sf d}\alpha&:=0,
\endaligned
\end{equation}
for the constant integer $\sf c$ as above.

Although slightly far from our main purpose in this section, but it may be worth to notice that the above structure equations, achieved as the output of the Cartan's algorithm, also reveal the so-called {\it structure constants} of the Lie algebra $\frak{aut}_{CR}(M({\bf a,b}))$ of infinitesimal CR-automorphisms associated to $M({\bf a,b})$. More details about the structure constants of a Lie algebra are presented in subsection \ref{structure constant}, below.

\begin{Proposition}
The Lie algebra $\frak{aut}_{CR}(M({\bf a,b}))$ of infinitesimal CR-automorphisms of $M({\bf a,b})$ is a $7$-dimensional algebra with the basis dual to the coframe $\{\mu,\sigma,\overline\sigma,\rho,\zeta,\overline\zeta,\alpha\}$ and with the structure constants as those visible in \thetag{\ref{struc-eq-prolonged}}.
\end{Proposition}

\subsection{Reformation of the coframe}
\label{reformation}

The constant invariants appeared in the latest structure equations \thetag{\ref{struc-eq-prolonged}} give us some {\it sufficient} conditions for the biholomorphic equivalence of two elements of the class $\mathcal C(1,4)$. But, to find a collection of necessary and sufficient invariants of this problem, one should reform the lifted coframe $(\mu,\sigma,\overline\sigma,\rho,\zeta,\overline\zeta)$ into a more adaptive one ({\it cf.} \cite[Chapter 8]{Olver-1995}). To this aim, first we need to divide the procedure into two distinct cases ${\bf b}\neq 0$ and ${\bf b}=0$.

 \subsubsection{The case ${\bf b}\neq 0$}

In this case and for the appropriate reformation of the current coframe, let us apply the substitution $\mu\mapsto {\sf c}'\,\mu^{\sf new}$ for ${\sf c}':=\frac{3}{2}{\sf c}$ that changes only the first structure equation of \thetag{\ref{struc-eq-prolonged}}:
\begin{equation}
\aligned
\label{struc-eq-reform-1}
{\sf d}\mu^{\sf new}&:=4\,\alpha\wedge\mu^{\sf new}+
\frac{\bf a}{\sf c'}\,\sigma\wedge\zeta+{\textstyle \frac{2}{3}}\,\sigma\wedge\overline\zeta+{\textstyle\frac{2}{3}}\,\overline\sigma\wedge\zeta+\frac{\overline{\bf a}}{\sf c'}\,\overline\sigma\wedge\overline\zeta,
\\
{\sf d}\sigma&:=3\,\alpha\wedge\sigma+\rho\wedge\zeta,
\\
{\sf d}\rho&:=2\,\alpha\wedge\rho+i\,\zeta\wedge\overline\zeta,
\\
{\sf d}\zeta&:=\alpha\wedge\zeta,
\\
{\sf d}\alpha&:=0,
\endaligned
\end{equation}
for which $\sf c'$ is clearly either $\bf b$ or $-\bf b$. But still, it is possible to have more reformation on the coframe. Let us change the three 1-forms:
\begin{equation}
\label{coframe-new}
\sigma\mapsto \frac{\sf c'}{\bf a}\,\sigma^{\sf new}, \ \ \ \rho\mapsto\frac{\sf c'}{\bf a}\,\rho^{\sf new}, \ \ \ \overline\zeta\mapsto\frac{\bf a}{\sf c'}\,\overline\zeta^{\sf new},
\end{equation}
and keep the remaining ones $\mu^{\sf new},\overline\sigma, \zeta,\alpha$ as before.
Then, the structure equations \thetag{\ref{struc-eq-reform-1}} convert into the fully-presented form:
\begin{equation}
\aligned
\label{struc-eq-reform-2}
{\sf d}\mu^{\sf new}&:=4\,\alpha\wedge\mu^{\sf new}+
\sigma^{\sf new}\wedge\zeta+{\textstyle \frac{2}{3}}\,\sigma^{\sf new}\wedge\overline\zeta^{\sf new}+{\textstyle\frac{2}{3}}\,\overline\sigma\wedge\zeta+\frac{{\bf a}\overline{\bf a}}{{\bf b}^2}\,\overline\sigma\wedge\overline\zeta^{\sf new},
\\
{\sf d}\sigma^{\sf new}&:=3\,\alpha\wedge\sigma^{\sf new}+\rho^{\sf new}\wedge\zeta,
\\
{\sf d}\overline\sigma&:=3\,\alpha\wedge\overline\sigma+\rho^{\sf new}\wedge\overline\zeta^{\sf new},
\\
{\sf d}\rho^{\sf new}&:=2\,\alpha\wedge\rho^{\sf new}+i\frac{{\bf a}^2}{{\bf b}^2}\,\zeta\wedge\overline\zeta^{\sf new},
\\
{\sf d}\zeta&:=\alpha\wedge\zeta,
\\
{\sf d}\overline\zeta^{\sf new}&:=\alpha\wedge\overline\zeta^{\sf new},
\\
{\sf d}\alpha&:=0.
\endaligned
\end{equation}
 Now one observes two invariants amongst the structure equations, namely the real invariant $\frak R:=\frac{{\bf a}\overline{\bf a}}{{\bf b}^2}$ and the complex one $\frak C:=\frac{{\bf a}^2}{{\bf b}^2}$. According to the simple facts in the theory of complex variables, one easily verifies that satisfying the latter invariance $\frak C$ guarantees it for the former one $\frak R$.

 Nevertheless, this fact should not deceive one to take $\frak C$ definitely as the single essential invariant of the problem. In fact, after substituting the initial 1-forms as \thetag{\ref{coframe-new}},  $\rho^{\sf new}$ is no longer a real 1-form. More precisely, here we have $\overline\rho^{\sf new}=\frac{\overline{\bf a}}{{\sf c}'}\rho=\frac{\overline{\bf a}}{\bf a}\rho^{\sf new}$ and hence $\overline\rho^{\sf new}$ is a constant multiple of $\rho^{\sf new}$. Consequently, two collections $\Theta:=\{\mu^{\sf new},\sigma^{\sf new},\overline\sigma,\rho^{\sf new},\zeta,\overline\zeta^{\sf new},\alpha\}$ and $\Omega:=\{\mu^{\sf new},\sigma^{\sf new},\overline\sigma,\overline\rho^{\sf new},\zeta,\overline\zeta^{\sf new},\alpha\}$ can be considered as the two very little different coframes on the prolonged space. On the other hand, according to the general principles behind the theory of the equivalence of coframes ({\it cf.} \cite[pp. 257-8]{Olver-1995}), one is not imposed to compute each structure equation in terms of the wedge product of a single coframe. Thus, let us replace in \thetag{\ref{struc-eq-reform-2}} the structure equation of ${\sf d} \rho^{\sf new}$ by ${\sf d}\overline\rho^{\sf new}$ expressed in terms of the 1-forms in $\Omega$:
\begin{equation}
\aligned
\label{struc-eq-reform-3}
{\sf d}\mu^{\sf new}&:=4\,\alpha\wedge\mu^{\sf new}+
\sigma^{\sf new}\wedge\zeta+{\textstyle \frac{2}{3}}\,\sigma^{\sf new}\wedge\overline\zeta^{\sf new}+{\textstyle\frac{2}{3}}\,\overline\sigma\wedge\zeta+\frak R\,\overline\sigma\wedge\overline\zeta^{\sf new},
\\
{\sf d}\sigma^{\sf new}&:=3\,\alpha\wedge\sigma^{\sf new}+\rho^{\sf new}\wedge\zeta,
\\
{\sf d}\overline\sigma&:=3\,\alpha\wedge\overline\sigma+\rho^{\sf new}\wedge\overline\zeta^{\sf new},
\\
{\sf d}\overline\rho^{\sf new}&:=2\,\alpha\wedge\overline\rho^{\sf new}+i\frak R\,\zeta\wedge\overline\zeta^{\sf new},
\\
{\sf d}\zeta&:=\alpha\wedge\zeta,
\\
{\sf d}\overline\zeta^{\sf new}&:=\alpha\wedge\overline\zeta^{\sf new},
\\
{\sf d}\alpha&:=0.
\endaligned
\end{equation}
As one observes, these final structure equations admits just the essential invariant $\frak R$. Then we can conclude the procedure by the following result;

\begin{Theorem}
\label{Th-approach-1}
In the case that ${\bf b}\neq 0$, the single essential invariant of the biholomorphic equivalence problem for the CR-manifolds belonging to the class $\mathcal C(1,4)$ comprising 6-dimensional Beloshapka's CR-models $M({\bf a,b})\subset\mathbb C^{1+4}$, represented as the graph of four polynomial functions:
\begin{equation*}
\aligned
w_1-\overline w_1&=2i\,z\overline z,
\\
w_2-\overline w_2&=2i\,(z^2\overline z+z\overline z^2),
\\
w_3-\overline w_3&=2\,(z^2\overline z-z\overline z^2),
\\
w_4-\overline w_4&=2i\,({\bf a}\,z^3\overline z+\overline{\bf a}\,z\overline z^3)+2i{\bf b}\,z^2\overline z^2,
\endaligned
\end{equation*}
is:
\[\boxed{
\frak R=\frac{{\bf a}\overline{\bf a}}{{\bf b}^2}.}
\]
\end{Theorem}

\begin{Remark}
Comparing to \cite{Beloshapka, Mamai}, one readily verifies that the achieved essential invariant $\frak R$ is exactly as that computed by Beloshapka and next by Mamai by means of the employed techniques from invariant theory.
\end{Remark}

\subsubsection{The case ${\bf b}=0$}

In this case, the latest structure equation \thetag{\ref{struc-eq-prolonged}} converts into the simpler form:
\begin{equation*}
\aligned
{\sf d}\mu&:=4\,\alpha\wedge\mu+
{\bf a}\,\sigma\wedge\zeta+\overline{\bf a}\,\overline\sigma\wedge\overline\zeta,
\\
{\sf d}\sigma&:=3\,\alpha\wedge\sigma+\rho\wedge\zeta,
\\
{\sf d}\rho&:=2\,\alpha\wedge\rho+i\,\zeta\wedge\overline\zeta,
\\
{\sf d}\zeta&:=\alpha\wedge\zeta,
\\
{\sf d}\alpha&:=0.
\endaligned
\end{equation*}
Now, by the substitutions $\sigma\mapsto\frac{1}{\bf a}\,\sigma^{\sf new}$, $\overline\sigma\mapsto\frac{1}{\overline{\bf a}}\,\overline\sigma^{\sf new}$ and $\rho\mapsto\frac{1}{\bf a}\,\rho^{\sf new}$ one receives the structure equations:
\begin{equation*}
\aligned
{\sf d}\mu&:=4\,\alpha\wedge\mu+\sigma^{\sf new}\wedge\zeta+\overline\sigma^{\sf new}\wedge\overline\zeta,
\\
{\sf d}\sigma^{\sf new}&:=3\,\alpha\wedge\sigma^{\sf new}+\rho^{\sf new}\wedge\zeta,
\\
{\sf d}\overline\sigma^{\sf new}&:=3\,\alpha\wedge\overline\sigma^{\sf new}+\overline\rho^{\sf new}\wedge\zeta,
\\
{\sf d}\rho&:=2\,\alpha\wedge\rho+i\,\zeta\wedge\overline\zeta,
\\
{\sf d}\zeta&:=\alpha\wedge\zeta,
\\
{\sf d}\alpha&:=0.
\endaligned
\end{equation*}
Here we encounter no any essential invariant among the structure equations. This means that:

\begin{Proposition}
\label{prop-b=0}
All CR-models of the form $M({\bf a},0), {\bf a}\neq 0$ are biholomorphically equivalent.
\end{Proposition}

\section{Second approach: Lie symmetry}
\label{approach-2}

Having in hand the structure of symmetry Lie algebras of the under consideration CR-models, namely their Lie algebras of infinitesimal CR-automorphisms in the terminology of CR-geometry, can help us to bypass a major part of the Cartan's algorithm to construct the associated structure equations. The key result of this approach is the following proposition of Beloshapka ({\it see} also \thetag{\ref{g-int}} and the paragraph containing it):

\begin{Proposition}
(see \cite[Proposition 3]{Beloshapka2004}). Each universal CR-model $M$ is isomorphic to the compact Lie group ${\sf Aut}_-(M)$ associated to the negative part, namely Levi-Tanaka algebra, $\frak{aut}_-(M)$ of $\frak{aut}_{CR}(M)$.
\end{Proposition}

This proposition suggests oneself to take into consideration an alternative equivalence problem instead of the original problem of equivalence between universal models. In fact by this proposition, the equivalence problem between two CR-models $M$ and $M'$ is manifesting itself as the equivalence problem of the two associated compact Lie groups ${\sf Aut}_-(M)$ and ${\sf Aut}_-(M')$. In this section, the idea is to investigate the equivalence problem of these associated Lie groups instead of the original equivalence problem for the models.

\subsection{Structure constants}
\label{structure constant}
Anyway, a natural question arises at this juncture is about the advantage of considering the equivalence problem between the Lie groups ${\sf Aut}_-(M)$ instead of universal CR-models $M$. To answer this question and as we know, if
$\{{\sf v}_1,\ldots, {\sf v}_r\}$ is a basis of a real Lie algebra
$\frak g$ of left-invariant vector fields on
an $r$-dimensional Lie group $G$, enjoying the Lie brackets:
\[
\big[{\sf v}_i,{\sf v}_j\big]
=
\sum_{k=1}^r\,c^k_{ij}\,{\sf v}_k
\ \ \ \ \ \ \ \ \ \ \ \ \
{\scriptstyle{(1\,\leqslant\,i\,<\,j\,\leqslant\,r)}},
\]
with certain {\sl structure constants}
$c^\bullet_{\bullet\bullet}$,
and if $\alpha^1,\ldots,\alpha^r$
is the dual Maurer-Cartan basis of left-invariant
1-forms, then their structure equations:
\begin{equation}
\label{struc-cons-general}
d\alpha^k
=
-\sum_{1\leqslant i<j\leqslant r}\,
c^k_{ij}\,\alpha^i\wedge\alpha^j
\ \ \ \ \ \ \ \ \ \ \ \ \
{\scriptstyle{(k\,=\,1\,\cdots\,r)}}
\end{equation}
have the same structure coefficients $c^\bullet_{\bullet\bullet}$
up to an overall minus sign.

{\it Thanks to this fact and in the case that the structure of the Lie algebra $\frak g$ is in hand, computing the structure equations of the associated compact Lie group $G$ is quite straightforward and does not necessitate one to perform the partly complicated parts of the Cartan's algorithm.}

\medskip
\noindent
{\bf Notations.} Let $M({\bf a,b})$ be a certain CR-model belonging to the class $\mathcal C(1,4)$. From now on, we denote by $\frak g({\bf a,b})$ and $G({\bf a,b})$ respectively the the Levi-Tanaka subalgebra $\frak{aut}_-(M({\bf a,b}))$ of $\frak{aut}_{CR}(M({\bf a,b}))$ and its associated Lie group ${\sf Aut}_-(M({\bf a,b}))$.

\subsection{Computing Lie algebras $\frak g({\bf a,b})$}

 Although this current idea of employing the already mentioned symmetry Lie algebras $\frak g({\bf a,b})$ seems, at the first look, simpler than that of performing the Cartan's algorithm as the preceding section, one should notice that here the major difficulty is to compute the basis elements of these Lie algebras. Roughly speaking, it seems that the complicated parts of the Cartan's algorithm manifest themselves in the computation of these symmetry algebras. As is mentioned in subsection \ref{preliminary-infinitesimal}, computing such algebras by means of the classical method of constructing and solving the so-called determining {\sc pde} systems is quite expensive, but recently, we have designed in \cite{SHAM} a powerful algorithm for computing such algebras which relies just on some simple techniques of linear algebra and some effective tools from the modern concept of comprehensive Gr\"obner systems. The implementation of this algorithm in the {\sc Maple} program is able compute $\frak g({\bf a,b})$ in just a few seconds. Before presenting the outputs, let us consider the following result which can simplify the subsequent computations:

\begin{Lemma}
\label{a-r}
Each 6-dimensional CR-model $M({\bf a,b})$ represented as the graph of four defining polynomials:
\begin{equation}
\label{model-2}
\aligned
w_1-\overline w_1&=2i\,z\overline z,
\\
w_2-\overline w_2&=2i\,(z^2\overline z+z\overline z^2),
\\
w_3-\overline w_3&=2\,(z^2\overline z-z\overline z^2),
\\
w_4-\overline w_4&=2i\,({\bf a}\,z^3\overline z+\overline{\bf a}\,z\overline z^3)+2i{\bf b}\,z^2\overline z^2,
\endaligned
\end{equation}
is equivalent to the CR-model $M({\sf r},{\bf b})$ where ${\sf r}:=({\bf a}\overline{\bf a})^{\frac{1}{2}}$.
\end{Lemma}

\proof
Consider the nonzero complex integer $\bf a$ in polar coordinates as ${\bf a}:={\sf r}e^{i\theta}$. Then the simple holomorphic transformation:
\[
z\mapsto e^{-i\frac{\theta}{2}}z
\]
converts the fourth defining polynomial of \thetag{\ref{model-2}} into the form:
\[
v_4=2i\,({\sf r}\,z^3\overline z+{\sf r}\,z\overline z^3)+2i{\bf b}\,z^2\overline z^2.
\]
In this case, one can transform also the first three equations of \thetag{\ref{model-2}} into their original forms by using some simple holomorphic transformations of the complex variables $w_1,w_2$ and $w_3$.
\endproof

Although we could use the above result to simplify the computations in section \ref{Cartan} but we proceed those computations without taking this very simplifying result into account since we wanted to show, as well, the effectiveness of the Cartan's algorithm.

 Lemma \ref{a-r} permits us to consider just the equivalence problem of the specific 6-dimensional CR-models $M({\sf r},{\bf b}), 0\neq {\sf r},{\bf b}\in\mathbb R$. For this aim, first we need the structure of the Lie algebra $\frak g({\sf r},{\bf b})$. According to the output\footnote{Actually we have computed this algebra for the general case $\frak g({\bf a,b})$ in which one can find the output in \cite{Maple}.} of the algorithm designed in \cite{SHAM}, $\frak g({{\sf r},{\bf b}})$ is a 6-dimensional real algebra generated by the vector fields ${\sf L},\widetilde{\sf L}, {\sf T}, {\sf S}, {\widetilde{\sf S}}, {\sf U}$, with the following table of Lie brackets:
\medskip
\begin{center}
\footnotesize
\begin{tabular} [t] { c | c | c | c | c | c | c }
& ${\sf L}$ & ${\widetilde{\sf L}}$ & ${\sf T}$ & ${\sf S}$ & $\widetilde{{\sf S}}$
& ${\sf U}$
\\
\hline ${\sf L}$ & $0$ & $-{\sf r}^2\,{\sf T}$ & $-2{\sf r}\,\widetilde{{\sf S}}$ & $0$ & $(-{\sf r}{\bf b}+\frac{3}{2}{\sf r}^2)\,{\sf U}$ & $0$
\\
\hline
$\widetilde{{\sf L}}$ & $*$ & $0$ & $-2{\sf r}\,{\sf S}$ & $(-{\sf r}{\bf b}-\frac{3}{2}{\sf r}^2)\,{\sf U}$ & $0$ & $0$
\\
\hline
${\sf T}$ & $*$ & $*$ & $0$ & $0$ & $0$ & $0$
\\
\hline
${\sf S}$ & $*$ & $*$ & $*$ & $0$ & $0$ & $0$
\\
\hline
$\widetilde{{\sf S}}$ & $*$ & $*$ & $*$ & $*$ & $0$ & $0$
\\
\hline
${\sf U}$ & $*$ & $*$ & $*$ & $*$ & $*$ & $0$ %
\end{tabular}
\end{center}
This Lie algebra is graded, in the sense of Tanaka:
\[
\frak g({\sf r},{\bf b}):=\frak g_{-4}\oplus\frak g_{-3}\oplus\frak g_{-2}\oplus\frak g_{-1},
\]
with $\frak g_{-1}={\rm Span}_{\mathbb R}\langle {\sf L, \widetilde{L}}\rangle$, with $\frak g_{-2}={\rm Span}_{\mathbb R}\langle {\sf T}\rangle$, with $\frak g_{-3}={\rm Span}_{\mathbb R}\langle {\sf S, \widetilde{S}}\rangle$ and with $\frak g_{-4}={\rm Span}_{\mathbb R}\langle {\sf U}\rangle$.
Here, we do not need the explicit expressions of these vector fields though one can find them in \cite{Maple}.

After providing the structure of the Lie algebra $\frak g({\sf r},{\bf b})$, now we are ready to plainly compute the structure equations associated to $G({\sf r}, {\bf b})$. Consider
\[
(\mu'',\sigma'',\widetilde{\sigma}'',\rho'',\zeta'',\widetilde{\zeta}'') \,\, \textrm{ as the dual coframe of}\,\, ({\sf U}, {\sf S},  {\widetilde{\sf S}},  {\sf T}, {\sf L},\widetilde{\sf L}).
 \]
Then according to \thetag{\ref{struc-cons-general}}, the structure equations of $G({\sf r},{\bf b})$ are:

\begin{equation}
\label{struc-approach2-1}
\aligned
{\sf d}\mu''&=({\sf r}{\bf b}-\textstyle\frac{3}{2}\,{\sf r}^2)\,\zeta''\wedge\widetilde{\sigma}''+({\sf r}{\bf b}+\frac{3}{2}\,{\sf r}^2)\,\widetilde{\zeta}''\wedge\sigma'',
\\
{\sf d}\sigma''&=2{\sf r}\,\widetilde{\zeta}''\wedge\rho'',
\\
{\sf d}\widetilde{\sigma}''&=2{\sf r}\,\zeta''\wedge\rho'',
\\
{\sf d}\rho''&={\sf r}^2\,\zeta''\wedge\widetilde{\zeta}'',
\\
{\sf d}\zeta''&=0,
\\
{\sf d}\widetilde{\zeta}''&=0.
\endaligned
\end{equation}
At this moment, proceeding along the same lines as subsection \ref{reformation} by applying the simple substitutions $\zeta''\mapsto \frac{1}{{\sf r}}\zeta'$ and $\widetilde{\zeta}''\mapsto\frac{1}{{\sf r}}\widetilde{\zeta}'$ and afterwards the substitutions:
\[
\aligned
\zeta'&\mapsto \zeta+\widetilde{\zeta}, \ \ \ \ \ \sigma''\mapsto \sigma+\widetilde{\sigma},
\\
\widetilde{\zeta}'&\mapsto\zeta-\widetilde{\zeta}, \ \ \ \ \ \widetilde{\sigma}''\mapsto -\sigma+\widetilde{\sigma},
\endaligned
\]
with dropping all $''$ from the notations, one receives the
structure equations \thetag{\ref{struc-approach2-1}} as:
\begin{equation}
\label{struc-approach2-2}
\aligned
{\sf d}\mu&=3{\sf r}\,\zeta\wedge\sigma+2{\bf b}\,\zeta\wedge\widetilde{\sigma}-2{\bf b}\,\widetilde{\zeta}\wedge\sigma-3{\sf r}\,\widetilde{\zeta}\wedge\widetilde{\sigma},
\\
{\sf d}\sigma&=-2\,\widetilde{\zeta}\wedge\rho,
\\
{\sf d}\widetilde{\sigma}&=2\,\zeta\wedge\rho,
\\
{\sf d}\rho&=-2\,\zeta\wedge\widetilde{\zeta},
\\
{\sf d}\zeta&=0,
\\
{\sf d}\widetilde{\zeta}&=0.
\endaligned
\end{equation}
Here, one observes the analogous of these structure equations to the output \thetag{\ref{struc-eq-prolonged}} of the Cartan's algorithm without the Maurer-Cartan form $\alpha$. Thanks to Proposition \ref{prop-b=0} and to save the space, let us continue with the assumption that ${\bf b}\neq 0$. Then, applying the substitution $\mu\mapsto 2{\bf b}\,\mu^{\sf new}$ and subsequently applying:
\[
\widetilde{\sigma}\mapsto\frac{2{\bf b}}{3\sf r}\,\widetilde{\sigma}^{\sf new}, \ \ \ \ \zeta\mapsto\frac{3\sf r}{2\bf b}\,\zeta^{\sf new}
\]
and finally using the structure equation of $\rho^{\sf new}:=\frac{3\sf r}{2\bf b}\rho$ instead of $d\rho$ in the fourth equation bring the following form for the structure equations:
\begin{equation}
\label{struc-approach2-3}
\aligned
{\sf d}\mu^{\sf new}&=\textstyle\frac{9{\sf r}^2}{4{\bf b}^2}\,\zeta^{\sf new}\wedge\sigma+\zeta^{\sf new}\wedge\widetilde{\sigma}^{\sf new}-\widetilde{\zeta}\wedge\sigma-\widetilde{\zeta}\wedge\widetilde{\sigma}^{\sf new},
\\
{\sf d}\sigma&=-2\,\widetilde{\zeta}\wedge\rho,
\\
{\sf d}\widetilde{\sigma}^{\sf new}&=\textstyle\frac{9{\sf r}^2}{2{\bf b}^2}\,\zeta^{\sf new}\wedge\rho,
\\
{\sf d}\rho^{\sf new}&=-2\,\zeta\wedge\widetilde{\zeta},
\\
{\sf d}\zeta&=0,
\\
{\sf d}\widetilde{\zeta}&=0,
\endaligned
\end{equation}
that leads us to the following result;

\begin{Proposition}
The single essential invariant for the equivalence problem of two 6-dimensional CR-models of the form $M({\sf r},{\bf b})$ for two real nonzero integers ${\sf r}$ and ${\bf b}$ is:
\[
\frak I:=\frac{{\sf r}^2}{{\bf b}^2}.
\]
\end{Proposition}

Now according to Lemma \ref{a-r}, each CR-model $M({\sf r}, {\bf b}), {\sf r},{\bf b}\in\mathbb R$ is biholomorphic equivalent to any CR-model $M({\bf a,b}), {\bf a}\in\mathbb C, {\bf b}\in\mathbb R$ when ${\sf r}=({\bf a}\overline{\bf a})^{\frac{1}{2}}$. Combining this fact to the above result gives the single invariant of our equivalence problem precisely as that of Theorem \ref{Th-approach-1};

\begin{Theorem}
\label{Th-approach-2}
The single essential invariant of the biholomorphic equivalence problem for the CR-manifolds belonging to the class $\mathcal C(1,4)$ comprising 6-dimensional Beloshapka's CR-models $M({\bf a,b})\subset\mathbb C^{1+4}$, for two nonzero integers ${\bf a}\in\mathbb C$ and ${\bf b}\in\mathbb R$, defined as the graph of four polynomial functions:
\begin{equation*}
\aligned
w_1-\overline w_1&=2i\,z\overline z,
\\
w_2-\overline w_2&=2i\,(z^2\overline z+z\overline z^2),
\\
w_3-\overline w_3&=2\,(z^2\overline z-z\overline z^2),
\\
w_4-\overline w_4&=2i\,({\bf a}\,z^3\overline z+\overline{\bf a}\,z\overline z^3)+2i{\bf b}\,z^2\overline z^2,
\endaligned
\end{equation*}
is:
\[\boxed{
\frak R=\frac{{\bf a}\overline{\bf a}}{{\bf b}^2}.}
\]
\end{Theorem}

\subsection*{Acknowledgment} The author would like to express his sincere thanks to Valerii Beloshapka and Jo\"el Merker for their helpful remarks, encouragements,
and suggestions during the preparation of this paper. He also would like to thank Benyamin Alizadeh for providing him the required Lie algebra $\frak g({\bf a,b})$ as the output of the implementation of the algorithm designed in \cite{SHAM}. The author also would like to thank Shahrekord University for its financial support.

\end{document}